\documentclass[12pt]{amsart}
\usepackage{amsfonts}
\usepackage{amsmath}

\theoremstyle{plain}
\newtheorem{X}{X}[section]
\newtheorem{theorem}[X]{Theorem}

\newtheorem{lemma}[X]{Lemma}
\newtheorem{corollary}[X]{Corollary}

\theoremstyle{definition}

\theoremstyle{remark}

\newtheorem{remark}[X]{Remark}

\numberwithin{equation}{section}
\numberwithin{equation}{section}

\newcommand\Aut{\operatorname{Aut}}
\newcommand\Gal{\operatorname{Gal}}
\newcommand\Sh{\operatorname{Sh}}
\newcommand\Spec{\operatorname{Spec}}
\def\ab{{ \text{ab} }}

\begin{document}
\title{Descent for Shimura Varieties}
\author{J.S. Milne}
\address{Department of Mathematics\\
University of Michigan\\
Ann Arbor, MI 48109}
\email{jmilne@umich.edu; http://www.math.lsa.umich.edu/$\sim$jmilne/}
\date{September 22, 1998}
\thanks{A small part of this research was supported by the National Science
Foundation.}

\begin{abstract}
This note proves that the descent maps provided by Langlands's Conjugacy
Conjecture \emph{do} satisfy the continuity condition necessary for them to
be effective. Hence the conjecture does imply the existence of canonical
models.
\end{abstract}

\maketitle

In his Corvallis article (1979, \S 6), Langlands stated a conjecture that
identifies the conjugate of a Shimura variety by an automorphism of $\mathbb{%
C}$ with the Shimura variety defined by different data, and he sketched a
proof that his conjecture implies the existence of canonical models.
However, as J\"{o}rg Wildeshaus and others have pointed out to me, it is not
obvious that the descent maps defined by Langlands satisfy the continuity
condition necessary for the descent to be effective. In this note, I prove
that they do satisfy this condition, and hence that Langlands's conjecture 
\emph{does\/} imply the existence of canonical models --- this is our only
proof of the existence of these models for a general Shimura variety. The
proof is quite short and elementary. I give it in Section 2 after reviewing
some generalities on the descent of varieties in Section~1.

\subsubsection{Notations and Conventions}

A variety over a field $k$ is a geometrically reduced scheme of finite type
over $\Spec k$ (\emph{not} necessarily irreducible). For a variety $V$ over
a field $k$ and a homomorphism $\sigma \colon k\rightarrow k^{\prime }$, $%
\sigma V$ is the variety over $k^{\prime }$ obtained by base change. The
ring of finite ad\`{e}les for $\mathbb{Q}$ is denoted by $\mathbb{A}_{f}{}$.

\section{Descent of Varieties.}

In this section, $\Omega$ is an algebraically closed field of characteristic
zero. For a field $L\subset\Omega$, $A(\Omega/L)$ denotes the group of
automorphisms of $\Omega$ fixing the elements of $L$.

Let $V$ be a variety over $\Omega$, and let $k$ be a subfield of $\Omega$. A
family $(f_{\sigma})_{\sigma\in A(\Omega/k)}$ of isomorphisms $f_{\sigma
}\colon\sigma V\rightarrow V$ will be called a \emph{descent system\/} if $%
f_{\sigma\tau}=f_{\sigma}\circ\sigma f_{\tau}$ for all $\sigma,\tau\in
A(\Omega/k)$. We say that a model $(V_{0},f\colon V_{0,\Omega}\rightarrow V)$
of $V$ over $k$ \emph{splits} $(f_{\sigma})$ if $f_{\sigma}=f\circ(\sigma
f)^{-1}$ for all $\sigma\in A(\Omega/k)$, and that a descent system is \emph{%
effective\/} if it is split by some model over $k$. The next theorem
restates results of Weil 1956.

\begin{theorem}
Assume that $\Omega $ has infinite transcendence degree over $k$. A descent
system $(f_{\sigma })$$_{\sigma \in A(\Omega /k)}$ on a quasiprojective
variety $V$ over $\Omega $ is effective if, for some subfield $L$ of $\Omega 
$ finitely generated over $k$, the descent system $(f_{\sigma })_{\sigma \in
A(\Omega /L)}$ is effective.
\end{theorem}

\begin{proof}
Let $k^{\prime}$ be the algebraic closure of $k$ in $L$ --- then $k^{\prime}$
is a finite extension of $k$ and $L$ is a regular extension of $k^{\prime}$.
Let $(V_{t^{\prime}},f^{\prime}\colon V_{t^{\prime},\Omega}\rightarrow V)$
be the model of $V$ over $L$ splitting $(f_{\sigma})_{\sigma\in A(\Omega/L)}$%
. Let $t\colon L\rightarrow k_{t}$ be a $k^{\prime}$-isomorphism from $L$
onto a subfield $k_{t}$ of $\Omega$ linearly disjoint from $L$ over $%
k^{\prime}$, and let $V_{t}=V_{t^{\prime}}\otimes_{L,t}k_{t}$. Zorn's Lemma
allows us to extend $t$ to an automorphism $\tau$ of $\Omega$ over $%
k^{\prime}$. The isomorphism 
\begin{equation*}
f_{t,t^{\prime}}\colon V_{t^{\prime},\Omega}\overset{f^{\prime}}{\rightarrow 
}V\overset{f_{\tau}^{-1}}{\rightarrow}\tau V\overset{(\tau f^{\prime})^{-1}}{%
\rightarrow}V_{t,\Omega}
\end{equation*}
is independent of the choice of $\tau$, is defined over $L\cdot k_{t}$, and
satisfies the hypothesis of Weil 1956, Theorem 6, which gives a model $(W,f)$
of $V$ over $k^{\prime}$ splitting $(f_{\sigma})_{\sigma\in A(\Omega
/k^{\prime})}$.

For $\sigma \in A(\Omega /k)$, $g_{\sigma }\overset{\text{df}}{=}f_{\sigma
}\circ \sigma f\colon \sigma W_{\Omega }\rightarrow V$ depends only on $%
\sigma |k^{\prime }$. For $k$-homomorphisms $\sigma ,\tau \colon k^{\prime
}\rightarrow \Omega $, define $f_{\tau ,\sigma }=g_{\tau }^{-1}\circ
g_{\sigma }\colon \sigma W\rightarrow \tau W$. Then $f_{\tau ,\sigma }$ is
defined over the Galois closure of $k^{\prime }$ in $\Omega $ and the family 
$(f_{\tau ,\sigma })$ satisfies the hypotheses of Weil 1956, Theorem 3,
which gives a model of $V$ over $k$ splitting $(f_{\sigma })_{\sigma \in
A(\Omega /k)}.$
\end{proof}

\begin{corollary}
Let $\Omega$, $k$, and $V$ be as in the theorem, and let $%
(f_{\sigma})_{\sigma\in A(\Omega/k)}$ be a descent system on $V$. If there
is a finite set $\Sigma$ of points in $V(\Omega)$ such that

\begin{enumerate}
\item  any automorphism of $V$ fixing all $P\in \Sigma $ is the identity
map, and

\item  there exists a subfield $L$ of $\Omega$ finitely generated over $k$
such that $f_{\sigma}(\sigma P)=P$ for all $P\in\Sigma$ and all $\sigma\in
A(\Omega/L)$,
\end{enumerate}

then $(f_{\sigma})_{\sigma\in A(\Omega/k)}$ is effective.
\end{corollary}

\begin{proof}
After possibly replacing the $L$ in (b) with a larger finitely generated
extension of $k$, we may suppose that $V$ has a model $(W,f)$ over $L$ for
which the points of $\Sigma$ are rational, i.e., such that for each $%
P\in\Sigma$, $P=f(P^{\prime})$ for some $P^{\prime}\in W(L)$. Now, for each $%
\sigma\in A(\Omega/L)$, $f_{\sigma}$ and $f\circ\sigma f^{-1}$ are both
isomorphisms $\sigma V\rightarrow V$ sending $\sigma P$ to $P$, and so
hypothesis (a) implies they are equal. Hence $(f_{\sigma})_{\sigma\in
A(\Omega/L)}$ is effective, and the theorem applies.
\end{proof}

\begin{remark}
\begin{enumerate}
\item  It is easy to construct noneffective descent systems. For example,
take $\Omega$ to be the algebraic closure of $k$, and let $V$ be a variety $%
k $. A one-cocycle $h\colon A(\Omega/k)\rightarrow\Aut(V_{\Omega})$ can be
regarded as a descent system --- identify $h_{\sigma}$ with a map $\sigma
V_{\Omega }=V_{\Omega}\rightarrow V_{\Omega}$. If $h$ is not continuous, for
example, if it is a homomorphism into $\Aut(V)$ whose kernel is not open,
then the descent system will not be effective.

\item  An example (Dieudonn\'{e} 1964, p 131) shows that the hypothesis that 
$V$ be quasiprojective in (1.1) is necessary unless the model $V_{0}$ is
allowed to be an algebraic space in the sense of M. Artin.

\item  Theorem 1.1 and its corollary replace Lemma 3.23 of Milne 1994, which
omits the continuity conditions.
\end{enumerate}
\end{remark}

\subsection{Application to moduli problems.}

Suppose we have a contravariant functor $\mathcal{M}$ from the category of
algebraic varieties over $\Omega $ to the category of sets, and equivalence
relations $\sim $ on each of the sets $\mathcal{M}(T)$ compatible with
morphisms. The pair $(\mathcal{M},\sim )$ is then called a \emph{moduli
problem} over $\Omega $. A $t\in T(\Omega )$ defines a map 
\begin{equation*}
m\mapsto m_{t}\overset{\text{df}}{=}t^{\ast }m\colon \mathcal{M}%
(T)\rightarrow \mathcal{M}(\Omega ).
\end{equation*}

\noindent A \emph{solution to the moduli problem }is a variety $V$ over $%
\Omega$ together with an isomorphism $\alpha\colon\mathcal{M}(\Omega
)/\!\!\sim\rightarrow V(\Omega)$ such that:

\begin{enumerate}
\item  for all varieties $T$ over $\Omega$ and all $m\in\mathcal{M}(T)$, the
map $t\mapsto\alpha(m_{t})\colon T(\Omega)\rightarrow V(\Omega)$ is regular
(i.e., defined by a morphism $T\rightarrow V$ of $\Omega$-varieties);

\item  for any variety $W$ over $\Omega$ and map $\beta\colon\mathcal{M}%
(\Omega)/\!\!\sim\rightarrow W(\Omega)$ satisfying the condition (a), the
map $P\mapsto\beta(\alpha^{-1}(P))\colon V(\Omega)\rightarrow W(\Omega)$ is
regular.
\end{enumerate}

\noindent Clearly, a solution to a moduli problem is unique up to a unique
isomorphism when it exists.

Let $(\mathcal{M},\sim )$ be a moduli problem over $\Omega $, and let $k$ be
a subfield $\Omega $. For $\sigma \in A(\Omega /k)$, define $^{\sigma }%
\mathcal{M}$ to be the functor sending an $\Omega $-variety $T$ to $\mathcal{%
M}(\sigma ^{-1}T)$. We say that $(\mathcal{M},\sim )$ is \emph{rational over}
$k$ if there is given a family $(g_{\sigma })_{\sigma \in A(\Omega /k)}$ of
isomorphisms $g_{\sigma }\colon ^{\sigma }\mathcal{M}\rightarrow \mathcal{M}$, compatible with $\sim$,
such that $g_{\sigma \tau }=g_{\sigma }\circ \sigma g_{\tau }$ for all $%
\sigma $, $\tau \in A(\Omega /k)$ --- the last equation means that $%
g_{\sigma \tau }(T)=g_{\sigma }(T)\circ g_{\tau }(\sigma ^{-1}T)$ for all
varieties $T$. Note that $^{\sigma }\mathcal{M}(\Omega )=\mathcal{M}(\Omega )
$, and that the rule $\sigma m=g_{\sigma }(m)$ defines an action of $%
A(\Omega /k)$ on $\mathcal{M}(\Omega )$. A \emph{solution to a moduli problem%
} $(\mathcal{M},\sim ,(g_{\sigma }))$ \emph{rational over} $k$ is a variety $%
V_{0}$ over $k$ together with an isomorphism $\alpha \colon \mathcal{M}%
(\Omega )/\!\!\sim \rightarrow V_{0}(\Omega )$ such that

\begin{enumerate}
\item  $(V_{0,\Omega},\alpha)$ is a solution to the moduli problem $(%
\mathcal{M},\sim)$ over $\Omega$, and

\item  $\alpha $ commutes with the actions of $A(\Omega /k)$ on $\mathcal{M}%
(\Omega )$ and $V_{0}(\Omega )$.
\end{enumerate}

\noindent Again, $(V_{0},\alpha)$ is uniquely determined up to a unique
isomorphism (over $k$) when it exists.

\begin{theorem}
Assume that $\Omega $ has infinite transcendence degree over $k$. Let $(%
\mathcal{M},\sim ,(g_{\sigma }))$ be a moduli problem rational over $k$ for
which $(\mathcal{M},\sim )$ has a solution $(V,\alpha )$ over $\Omega $.
Then $(\mathcal{M},\sim ,(g_{\sigma }))$ has a solution over $k$ if there
exists a finite subset $\Sigma \subset \mathcal{M}(\Omega )$ such that

\begin{enumerate}
\item  any automorphism of $V$ fixing $\alpha(P)$ for all $P\in\Sigma$ is
the identity map, and

\item  there exists a subfield $L$ of $\Omega $ finitely generated over $k$
such that $g_{\sigma }(P)\sim P$ for all $P\in \Sigma $ and all $\sigma \in
A(\Omega /L)$.
\end{enumerate}
\end{theorem}

\begin{proof}
The family $(g_{\sigma})$ defines a descent system on $V$, which Corollary
1.2 shows to be effective.
\end{proof}

\section{Descent of Shimura Varieties.}

In this section, all fields will be subfields of $\mathbb{C}$. For a
subfield $E$ of $\mathbb{C}$, $E^{\text{ab}}$ denotes the composite of all
the finite abelian extensions of $E$ in $\mathbb{C}$.

Let $(G,X)$ be a pair satisfying the axioms (2.1.1.1--2.1.1.3) of Deligne
1979 to define a Shimura variety, and let $\Sh(G,X)$ be the corresponding
Shimura variety over $\mathbb{C}$. We regard $\Sh(G,X)$ as a pro-variety
endowed with a continuous action of $G(\mathbb{A}_{f})$ --- in particular
(ibid. 2.7.1) this means that $\Sh(G,X)$ is a projective system of varieties 
$(\Sh_{K}(G,X))$ indexed by the compact open subgroups $K$ of $G(\mathbb{A}%
_{f})$. Let $[x,a]$ $=([x,a]_{K})_{K}$ denote the point in $\Sh(G,X)(\mathbb{%
C})$ defined by a pair $(x,a)\in X\times G(\mathbb{A}_{f})$, and let $E(G,X)$
be the reflex field of $(G,X)$. For a special point $x\in X$, let $%
E(x)\supset E(G,X)$ be the reflex field for $x$ and let 
\begin{equation*}
r_{x}\colon \Gal(E(x)^{\ab}/E(x))\rightarrow T(\mathbb{A}_{f})/T(\mathbb{Q}%
)^{-}
\end{equation*}
be the reciprocity map defined in Milne 1992, p164 (inverse to that in
Deligne 1979, 2.2.3). Here $T$ is a subtorus of $G$ such that $\text{Im}%
(h_{x})\subset T_{\mathbb{R}}$ and $T(\mathbb{Q})^{-}$ is the closure of $T(%
\mathbb{Q})$ in $T(\mathbb{A}{}_{f})$. A\emph{\ model\/} of $\Sh(G,X)$ over
a field $k$ is a pro-variety $S$ over $k$ endowed with an action of $G(%
\mathbb{A}_{f})$ and a $G(\mathbb{A}_{f})$-equivariant isomorphism $f\colon
S_{\mathbb{C}}\rightarrow \Sh(G,X)$. A model of $\Sh(G,X)$ over $E(G,X)$ is 
\emph{canonical } if, for each special point $x\in X$ and $a\in G(\mathbb{A}%
_{f})$, $[x,a]$ is rational over $E(x)^{\text{ab}}$ and $\sigma \in \Gal%
(E(x)^{\text{ab}}/E(x))$ acts on $[x,a]$ according\footnote{%
More precisely, the condition for $(S,f)$ to be canonical is the following:
if $P\in S(\mathbb{C})$ corresponds under $f$ to $[x,a]$, then $\sigma P$
corresponds under $f$ to $[x,r_{x}(\sigma )\cdot a].$} to the rule: 
\begin{equation*}
\sigma \lbrack x,a]=[x,r_{x}(\sigma )\cdot a].
\end{equation*}
Let $k$ be a field containing $E(G,X)$. A \emph{descent system} for $\Sh(G,X)
$ over $k$ is a family of isomorphisms 
\begin{equation*}
(f_{\sigma }\colon \sigma \Sh(G,X)\rightarrow \Sh(G,X))_{\sigma \in A(%
\mathbb{C}/k)}
\end{equation*}
such that,

\begin{enumerate}
\item  for all $\sigma,\tau\in A(\mathbb{C}/k)$, $f_{\sigma\tau}=f_{\sigma
}\circ\sigma f_{\tau}$, and

\item  for all $\sigma\in A(\mathbb{C}/k) $, $f_{\sigma}$ is equivariant for
the actions of $G(\mathbb{A}_{f})$ on $\Sh(G,X)$ and $\sigma\Sh(G,X)$.
\end{enumerate}

\noindent We say that a model $(S,f)$ of $\Sh(G,X)$ over $k$ \emph{splits }$%
(f_{\sigma})$ if $f_{\sigma}=f\circ\sigma f^{-1}$ for all $\sigma\in
A(\Omega/k)$, and that a descent system if \emph{effective} if it is split
by some model over $k$. A descent system $(f_{\sigma})$ for $\Sh(G,X)$ over $%
E(G,X)$ is \emph{canonical} if 
\begin{equation*}
f_{\sigma}(\sigma\lbrack x,a])=[x,r_{x}(\sigma|E(x)^{\ab})\cdot a]
\end{equation*}

\noindent whenever $x$ is a special point of $X$, $\sigma\in A(\mathbb{C}%
/E(x))$, and $a\in G(\mathbb{A}_{f})$.

\begin{remark}
\begin{enumerate}
\item  For a Shimura variety $\Sh(G,X)$, there exists at most one canonical
descent system for $\Sh(G,X)$ over $E(G,X)$. (Apply Deligne 1971, 5.1, 5.2.)

\item  Let $(S,f)$ be a model of $\Sh(G,X)$ over $E(G,X)$, and let $%
f_{\sigma }=f\circ(\sigma f)^{-1}$. Then $(f_{\sigma})_{\sigma\in A(\mathbb{C%
}/k)}$ is a descent system for $\Sh(G,X)$, and $(f_{\sigma})$ is canonical
if and only if $(S,f)$ is canonical.

\item  Suppose $\Sh(G,X)$ has a canonical descent system $%
(f_{\sigma})_{\sigma\in A(\mathbb{C}/E(G,X))}$; then $\Sh(G,X)$ has a
canonical model if and only if $(f_{\sigma})$ is effective. (Follows from
(a) and (b).)

\item  A descent system $(f_{\sigma})_{\sigma\in A(\mathbb{C}/k)}$ on $\Sh%
(G,X)$ defines for each compact open subgroup $K$ of $G(\mathbb{A}_{f})$ a
descent system $(f_{\sigma,K})_{\sigma\in A(\mathbb{C}/k)}$ on the variety $%
\Sh_{K}(G,X)$ (in the sense of \S1). If $(f_{\sigma})$ is effective, then so
also is $(f_{\sigma,K})$ for all $K$; conversely, if $(f_{\sigma,K})_{\sigma%
\in A(\mathbb{C}/k)}$ is effective (in the sense of \S1) for all
sufficiently small $K$, then $(f_{\sigma})_{\sigma\in A(\mathbb{C}/k)}$ is
effective (in the sense of this section).
\end{enumerate}
\end{remark}

\begin{lemma}
The automorphism group of the quotient of a bounded symmetric domain by a
neat arithmetic group is finite.
\end{lemma}

\begin{proof}
According to Mumford 1977, Proposition 4.2, such a quotient is an algebraic
variety of logarithmic general type, which implies that its automorphism
group is finite (Iitaka 1982, 11.12).

Alternatively, one sees easily that the automorphism group of the quotient
of a bounded symmetric domain $D$ by a neat arithmetic subgroup $\Gamma$ is $%
N/\Gamma$ where $N$ is the normalizer of $\Gamma$ in $\Aut(D)$. Now $N$ is
countable and closed (because $\Gamma$ is closed), and hence is discrete
(Baire category theorem). Because the quotient of $\Aut(D)$ by $\Gamma$ has
finite measure, this implies that $\Gamma$ has finite index in $N$. Cf.
Margulis 1991, II 6.3.
\end{proof}

\begin{theorem}
Every canonical descent system on a Shimura variety is effective.
\end{theorem}

\begin{proof}
Let $(f_{\sigma})_{\sigma\in A(\mathbb{C}/E(G,X))}$ be a canonical descent
system for the Shimura variety $\Sh(G,X)$. Let $K$ be a compact open
subgroup of $G(\mathbb{A}_{f})$, chosen so small that the connected
components of $\Sh_{K}(G,X)$ are quotients of bounded symmetric domains by 
\emph{neat} arithmetic groups. Let $x$ be a special point of $X$. According
to Deligne 1971, 5.2, the set $\Sigma=\{[x,a]_{K}\mid a\in G(\mathbb{A}%
_{f})\}$ is Zariski dense in $\Sh_{K}(G,X)$. Because the automorphism group
of $\Sh
_{K}(G,X)$ is finite, there is a finite subset $\Sigma_{f}$ of $\Sigma$ such
that any automorphism $\alpha$ of $\Sh_{K}(G,X)$ fixing each $P\in\Sigma_{f}$
is the identity map.

The rule 
\begin{equation*}
\sigma \ast \lbrack x,a]_{K}=[x,r_{x}(\sigma )\cdot a]_{K}
\end{equation*}
defines an action of $\Gal(E(x)^{\text{ab}}/E(x))$ on $\Sigma $ for which
the stabilizer of each point of $\Sigma $ is open. Therefore, there exists a
finite abelian extension $L$ of $E(x)$ such that $\sigma \ast P=P$ for all $%
P\in \Sigma _{f}$ and all $\sigma \in \Gal(E(x)^{\text{ab}}/L)$.

Now, because $(f_{\sigma})_{\sigma\in A(\Omega/E(G,X))}$ is canonical, $%
f_{\sigma,K}(\sigma P)=P$ for all $P\in\Sigma_{f}$ and all $\sigma\in A(%
\mathbb{C}{}/L)$, and we may apply Corollary 1.2 to conclude that $%
(f_{\sigma,K})_{\sigma\in A(\mathbb{C}/E(G,X))}$ is effective. As this holds
for all sufficiently small $K$, $(f_{\sigma})_{\sigma\in A(\mathbb{C}%
/E(G,X))}$ is effective.
\end{proof}

\begin{remark}
\begin{enumerate}
\item  If Langlands's Conjugacy Conjecture (Langlands 1979, p232, 233) is
true for a Shimura variety $\Sh(G,X)$, then $\Sh(G,X)$ has a canonical
descent system (ibid. \S 6; also Milne and Shih 1982, \S 7).

\item  Langlands's Conjugacy Conjecture is true for all Shimura varieties
(Milne 1983). Hence canonical models exist for all Shimura varieties.
\end{enumerate}
\end{remark}

Another proof, based on different ideas, that the descent maps given by Langlands's conjecture are effective can be found in Moonen 1998. (I thank the referee for this reference.)

\section*{References}

Deligne, P., Travaux de Shimura, in S\'{e}minaire Bourbaki, 23\`{e}me
ann\'{e}e (1970/71), Exp. No. 389, 123--165. Lecture Notes in Math., 244,
Springer, Berlin, 1971.

Deligne, P., Vari\'{e}t\'{e}s de Shimura: Interpr\'{e}tation modulaire, et
techniques de construction de mod\`{e}les canoniques, Proc. Symp. Pure Math.
33 Part 2, pp. 247--290, 1979.

Dieudonn\'{e}, J., Fondements de la G\'{e}om\'{e}trie Alg\'{e}brique
Moderne, Presse de l'Universit\'{e} de Montr\'{e}al, 1964.

Iitaka, S., Algebraic Geometry, Springer, Heidelberg, 1982.

Langlands, R., Automorphic representations, Shimura varieties, and motives,
Ein M\"{a}rchen, Proc. Symp. Pure Math. 33 Part 2, pp. 205--246, 1979.

Margulis, G.A., Discrete subgroups of semisimple Lie groups, Springer,
Heidelberg, 1991.

Milne, J.S., The action of an automorphism of $\mathbb{C}$ on a Shimura
variety and its special points, Prog. in Math., vol. 35, Birkh\"{a}user,
Boston, pp. 239--265, 1983.

Milne, J.S., The points on a Shimura variety modulo a prime of good
reduction, in The Zeta Function of Picard Modular Surfaces (Langlands and
Ramakrishnan, eds), Les Publications CRM, Montr\'{e}al, pp. 153--255, 1992.

Milne, J. S., Shimura varieties and motives, in Motives (Seattle,
WA, 1991), 447--523, Proc. Sympos. Pure Math., Part 2, Amer. Math. Soc.,
Providence, RI, 1994.

Milne, J.S. and Shih, K-y., Conjugates of Shimura varieties, in Hodge
Cycles, Motives, and Shimura Varieties, Lecture Notes in Math., vol. 900,
Springer, Heidelberg, pp. 280--356, 1982.

Moonen, B., Models of Shimura varieties in mixed characteristics, in Galois Representations in Arithmetic Geometry (A.J. Scholl and R.L Taylor, editors), Cambridge University Press, pp. 267--350, 1998.

Mumford, D., Hirzebruch's proportionality theorem in the non-compact case,
Invent. Math. 42, 239--272, 1977.

Weil, A., The field of definition of a variety, Amer. J. Math 78, pp.
509--524, 1956.

\end{document}